\documentclass{article}

\usepackage{arxiv}
\usepackage{amssymb}
\usepackage{amsmath}
\usepackage{empheq}
\usepackage{amscd}
\usepackage{graphicx}
\usepackage{graphics}
\usepackage[noadjust]{cite}
\usepackage{amsthm}
\usepackage{caption}
\usepackage{multirow}
\usepackage{array}
\usepackage{rotating}

\usepackage{indentfirst}
\usepackage{tikz}
\usepackage{calc}
\usepackage{epsfig}
\usepackage{natbib}
\usepackage[makeroom]{cancel}
\usepackage{multicol}
\usepackage{csquotes}
\usepackage{algorithm}
\usepackage{algorithmic}
\usepackage{enumitem}
\usepackage{hyperref}
\usepackage{url}
\hypersetup{colorlinks,linkcolor={blue},citecolor={blue},urlcolor={blue}}

\usepackage{timet}
\usepackage{graphicx}
\usepackage{multirow}
\usepackage{multicol}
\usepackage{lipsum}
\usepackage{timet}
\usepackage{epsfig}
\usepackage{amsmath}
\usepackage{amsfonts}
\usepackage{amssymb}
\usepackage{color}
\numberwithin{equation}{section}
\usepackage{caption}
\usepackage{float}
\usepackage{subcaption}
\usepackage{graphics}
\usepackage{xcolor,graphicx}
\usepackage{csquotes}
\usepackage{mathtools}
\usepackage{optidef}

\newcommand{\vertiii}[1]{{\left\vert\kern-0.25ex\left\vert\kern-0.25ex\left\vert #1 
    \right\vert\kern-0.25ex\right\vert\kern-0.25ex\right\vert}}

\def\tsc#1{\csdef{#1}{\textsc{\lowercase{#1}}\xspace}}
\tsc{WGM}
\tsc{QE}
\tsc{EP}
\tsc{PMS}
\tsc{BEC}
\tsc{DE}
\renewcommand{\u}[0]{\bold{u}}                        

\renewcommand{\b}[0]{\bold{b}}                        
\renewcommand{\d}[0]{\bold{d}}                        
\renewcommand{\P}[0]{\bold{P}}                        
\renewcommand{\next}[1]{#1_{(k)}}                         
\newcommand{\prev}[1]{#1_{(k-1)}}                       
\newcommand{\nextin}[1]{#1_{(l)}}                         
\newcommand{\previn}[1]{#1_{(l-1)}}                       

\newcommand{\m}[0]{\bold{m}}                        
\newcommand{\w}[0]{\omega}                          
\newcommand{\A}[0]{\bold{A}}                        
\newcommand{\del}[0]{\bold{\nabla}^2}               
\newcommand{\Diag}[0]{\text{Diag}}                  
\newcommand{\AL}[0]{\mathcal{L}}                    
\newcommand{\penaltyparb}[0]{\lambda}               
\newcommand{\penaltypard}[0]{\mu}                   
\newcommand{\dualb}[0]{\bold{v}}                    
\newcommand{\duald}[0]{\bold{w}}                    
\newcommand{\sdualb}[0]{\bold{b}}                   
\newcommand{\sduald}[0]{\bold{d}}                   
\newcommand{\reg}[0]{\mathcal{R}}                   
\newcommand{\aA}[0]{\bold{A}_{\epsilon}}            
\newcommand{\M}[0]{\bold{M}}                   
\newcommand{\I}[0]{\bold{I}}                   
\newcommand{\Nr}[0]{n_r}                         
\newcommand{\Ns}[0]{n_s}                       
\newcommand{\N}[0]{n}                          
\newcommand{\F}[0]{\bold{F}}                   
\renewcommand{\c}[0]{\bold{c}}                        
\newcommand{\iu}[0]{i}                   

\begin{document}
\title{Large-scale highly-accurate extended full waveform inversion using convergent Born series}

\author{\href{http://orcid.org/0000-0003-1805-1132}{\includegraphics[scale=0.06]{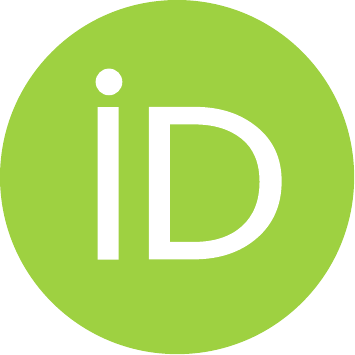}\hspace{1mm}Hossein S. Aghamiry} \\
  University Cote d'Azur - CNRS - IRD - OCA, Geoazur, Valbonne, France. 
  \texttt{aghamiry@geoazur.unice.fr}
\And
  \href{https://orcid.org/0000-0002-9879-2944}{\includegraphics[scale=0.06]{orcid.pdf}\hspace{1mm}Ali Gholami} \\
  Institute of Geophysics, University of Tehran, Tehran, Iran.
  \texttt{agholami@ut.ac.ir} \\ 
\And
 \href{https://orcid.org/0000-0002-9879-2944}{\includegraphics[scale=0.06]{orcid.pdf}\hspace{1mm}Kamal Aghazade} \\
  Institute of Geophysics, University of Tehran, Tehran, Iran.
  \texttt{aghazade.Kamal@ut.ac.ir} \\ 
\And
 \href{https://orcid.org/0000-0002-9879-2944}{\includegraphics[scale=0.06]{orcid.pdf}\hspace{1mm}Mahdi Sonbolestan} \\
  Institute of Geophysics, University of Tehran, Tehran, Iran.
  \texttt{momahdiso@yahoo.com} \\ 
  \And
\href{http://orcid.org/0000-0002-4981-4967}{\includegraphics[scale=0.06]{orcid.pdf}\hspace{1mm}St\'ephane Operto} \\ 
  University Cote d'Azur - CNRS - IRD - OCA, Geoazur, Valbonne, France. 
  \texttt{operto@geoazur.unice.fr} 
 }

\renewcommand{\shorttitle}{Extended full waveform inversion using convergent Born series, Aghamiry et al.}

\maketitle

\begin{abstract}
Full-waveform inversion (FWI) is a seismic imaging method that provides quantitative inference about subsurface properties with a wavelength-scale resolution. Its frequency-domain formulation is computationally efficient when processing only a few discrete frequencies. However, classical FWI, which is formulated on the reduced-parameter space, requires starting the inversion with a sufficiently-accurate initial model and low frequency to prevent being stuck in local minima due to cycle skipping. 
FWI with extended search space has been proposed to mitigate this issue. It contains two main steps: first, data-assimilated (DA) wavefields are computed by allowing for wave-equation errors to match the data at receivers closely. Then, subsurface parameters are estimated from these wavefields by minimizing the wave-equation errors.  
The DA wavefields are the least-squares solution of an overdetermined system gathering the wave and observation equations.
The numerical bandwidth of the resulting normal-equation system is two times that of the wave-equation system, which can be a limiting factor for 3D large-scale applications. Therefore, computing highly accurate DA wavefields at a reasonable computational cost is an issue in extended FWI. 
This issue is addressed here by rewriting the normal system such that its solution can be computed by solving the time-harmonic wave equation several times in sequence. Moreover, the computational burden of multi-right-hand side (RHS) simulations is mitigated with a sketching method. Finally, we solve the time-harmonic wave equation with the convergent Born series method, which conciliates accuracy and computational efficiency for problems whose size would be prohibitive for direct solvers. Application of the new extended FWI algorithm on the 2004 BP salt benchmark shows that it reconstructs at a 
reasonable cost subsurface models that are similar to those obtained with the classical extended FWI algorithm while being not limited by the size of the problem.
\end{abstract}

\section{Introduction}
Visualizing and quantifying the properties of a medium from sparse indirect measurements is the final goal of all the imaging methods in different application fields. 
Over the last decade, the so-called Full Waveform Inversion (FWI) has reached maturity gates and emerged as the leading-edge seismic imaging technique in exploration geophysics, along with the continuous development of high-performance computing and new acquisition technologies and geometries (dense node arrays, broadband sources, multi-component sensors) \citep{Virieux_2017_FWI}. This trend has contributed to disseminate FWI towards new fields of application such as earthquake seismology, microseismic relocation, time-lapse imaging for monitoring purposes, surface-wave interferometry, ambient noise imaging, distributed acoustic sensing (DAS) inversion, ground-penetrating radar (GPR), controlled-source electromagnetic inversion (CSEM), medical imaging and helioseismology among others.\\
FWI can be implemented in the time-space \cite{Tarantola_1984_ISR} and frequency-space \cite{Pratt_1998_GNF} domains, which have their own pros and cons \citep{Plessix_2017_CAT}. The latter is computationally efficient when the inversion can be limited to a few discrete frequencies, as is the case when dealing with data acquisition devices providing a broad angular illumination of the target \cite{Sirgue_2004_EWI}. Also, the frequency domain provides the most natural framework to implement multi-scale approaches by frequency continuation, which helps to mitigate the nonlinearity of the problem. Finally, the implementation of attenuation effects is straightforward and doesn't generate computational overhead in the frequency domain. \\ 
FWI is a partial differential equation (PDE)-constrained data-fitting problem for parameter estimation, where the constraint is the wave equation and the optimization variables gather the wavefields and the subsurface parameters, the latter being contained in the wave-equation coefficients \citep{Epanomeritakis_2008_NCG}.
Classical FWI is, however, formulated on the reduced parameter space after the projection of the wavefields in the data misfit function  \cite{Pratt_1998_GNF}. This variable projection leads to a highly nonlinear problem due to the oscillatory nature of waves. This nonlinearity manifests by the so-called cycle-skipping pathology, which traps the inversion in spurious minima when the starting model and the initial frequency don't predict the phase of the recorded arrivals with an error lower than half a period, a direct result of the Born approximation that is used to linearize the inverse problem. \\
Recently, different approaches have been proposed to bypass this shortage \cite[and the references therein]{Gholami_2022_EFW}. A broad category of methods extends the linear regime of FWI by enlarging the search space with non-physical degrees of freedom. Among these approaches, \cite{VanLeeuwen_2013_MLM} proposes the so-called frequency-domain wavefield-reconstruction inversion (WRI) method, which processes the wave equation as a soft constraint with a penalty method such that the simulated wavefields closely match the data at receivers, hence preventing cycle skipping. Then the subsurface parameters are updated by minimizing the wave-equation errors to push back the inversion toward the physics. Later, \citet{Aghamiry_2019_IWR} improved WRI by replacing the penalty method with an augmented Lagrangian method such that the penalty parameter can be kept fixed in iterations while the update of the Lagrange multipliers with gradient ascent steps guarantees that the constraint is satisfied at the convergence point with sufficient accuracy. This approach was called iteratively-refined (IR)-WRI. \\
In the (IR-)WRI framework, the wavefields are the least-squares solutions of an overdetermined system gathering the observation equation and the wave equation, the latter being weighted by the penalty parameter. These wavefields were referred to as data-assimilated (DA) wavefields by \citet{Aghamiry_2019_AEW} since they are computed with a feedback term to the observables.
%
%
%
%
%
%
%
%
The DA wavefields are the solution of a normal-equation system, which is denser, has a wider numerical bandwidth, and a poorer condition number than the original wave-equation system involved in classical FWI, where the wave-equation constraint is strictly satisfied at each iteration. These embarrassing specifications of the normal system can constitute limiting factors to tackle large-scale 3D problems whatever the discretized normal equation is solved with direct or iterative methods suitable for sparse linear systems. \\
In this study, we adapt in the frequency domain the new formulation of IR-WRI that has been recently proposed by \cite{Gholami_2022_EFW} to implement the normal equation in the time domain with explicit time-stepping schemes. This new formulation reformulates the normal equation such that the DA wavefields can be computed by repeated solutions of the wave equation and its adjoint. The main drawback of this approach is that the wave equation needs to be solved with an increased number of right-hand sides (RHSs) compared to the original formulation of WRI based on the direct solution of the normal equation. Typically, this number of RHSs scales to the number of sources + the number of receivers in the new formulation, while it scales to the number of sources only in the classical implementation. 
To address this issue, we use sketching methods to reduce the number of wave-equation solutions \cite{Aghazade_2021_RSS}. Such methods decrease the number of RHSs through a linear combination of sources/receivers with different weights at the price of extra artefacts in the reconstructed model and wavefields. Finally, we use the fast, highly-accurate and limited-memory convergent Born series (CBS) method \cite{Osnabrugge_2016_CBS} to further mitigate the computational burden of these simulations while being able to address large-scale simulations. Indeed, CBS is not the only possible choice and any kind of wave-equation solvers in the frequency or time domains can be interfaced with the proposed algorithm. 
In the time-domain case, a few monochromatic wavefields can be integrated on the fly in the loop over time steps  by discrete Fourier transform or phase sensitive detection to manage compact volume of data and implement the parameter-estimation subproblem in the frequency domain \cite{Nihei_2007_FRM,Sirgue_2008_FDW}. \\
This paper is mainly organized in a method, numerical example, discussion, and conclusion sections.
The method section begins with a review of the original frequency-domain IR-WRI formulation \cite{Aghamiry_2019_IWR} before considering its new formulation that bypasses the solution of the normal equation  \cite{Gholami_2022_EFW}. Then, we briefly review how the sketching method decreases the number of sources/receivers in the new implementation of IR-WRI. We finish the method section by describing the CBS method that is used to solve the wave equation and its adjoint. We continue with the numerical results section to assess the performance of the method against the challenging large-contrast 2004 BP salt model. We discuss further these results, some limitations of the proposed method, and some remedies before conclusions.
%
%
%
\section{Theory}
Full waveform inversion (FWI) with a general (convex) regularization term can be written as the following PDE constrained optimization \cite{Aghamiry_2019_IBC}:
\begin{mini} 
{\m,\u}{\reg(\m)}
{\label{eq:main}}{}
\addConstraint {\A(\m)\u}{=\b}
\addConstraint {\P\u}{=\d},
\end{mini}
where $\m \in \mathbb{R}^{\N\times 1}$ gathers the model parameters (the square of slowness), $\reg$ is a regularization function such as Tikhonov, total variation (TV) regularizers or their combination \cite{Aghamiry_2019_CRO}, or more general adaptive regularizers \cite{Aghamiry_2020_FWI},  $\A(\m)  \in \mathbb{C}^{\N\times \N}$ is the discretized Helmholtz operator,  $\u \in \mathbb{C}^{\N\times 1}$ is the state variable (wavefield), $\b \in \mathbb{C}^{\N\times 1}$ is the source term, $\P \in \mathbb{R}^{\Nr\times \N}$  is the observation operator that samples the simulated wavefield at the $\Nr$ receiver locations, and $\d \in \mathbb{C}^{\Nr\times 1}$ is the recorded wavefield at receiver locations (data). We assume a stationary-recording acquisition device, hence $\P$ is the same for all $\Ns$ sources.
The Helmholtz equation for angular frequency $\w$ reads
\begin{equation} \label{eq:PDE}
\del \u+ \w^2\Diag(\m)\u=\b,
\end{equation}
where $\del$ is the Laplace operator and $\Diag(\cdot)$ denotes a diagonal matrix with $\cdot$ on its main diagonal. We review the method for one source and one frequency. However, the algorithm easily generalizes to multisource and multifrequency configurations.

The augmented Lagrangian function associated with problem \eqref{eq:main} is given by \cite{Nocedal_2006_NO,Aghamiry_2019_IWR}
\begin{equation} \label{AL}
\AL(\m,\u,\dualb,\duald)= \reg(\m) + \langle\dualb,\b-\A(\m)\u\rangle + \langle\duald,\d-\P\u\rangle +\frac{\penaltyparb}{2} \|\b-\A(\m)\u\|_2^2 + \frac{\penaltypard}{2} \|\d -\P\u\|_2^2, 
\end{equation}
where $\langle\cdot,\cdot\rangle$ denotes inner product, $\dualb,\duald$ are the Lagrange multipliers, and $\penaltyparb,\penaltypard$ are the penalty parameters.
IR-WRI solves  \eqref{eq:main} via the following ADMM iteration \cite{Chen_1994_PDM,Boyd_2011_DOS,Aghamiry_2019_IWR}:
\begin{subequations}
\begin{align}
\next{\m} &= \arg\min_{\m} \mathcal{L}(\bold{\m},\prev{\u},\prev{\dualb},\prev{\duald}), \\
\next{\u} &= \arg\min_{\u} \mathcal{L}(\next{\m},\u,\prev{\dualb},\prev{\duald}),\\
\next{\dualb}& = \prev{\dualb} + \penaltyparb (\b -\A(\next{\m})\next{\u}),\\
\next{\duald}& = \prev{\duald} + \penaltypard (\d -\P\next{\u}).
\end{align}
\end{subequations}
We use a scaled form of the augmented Lagrangian method to simplify the ADMM iterations \citep[][Section 3.1.1]{Boyd_2011_DOS}. Introducing the scaled Lagrange multipliers $\next{\sdualb}= \next{\dualb}/\penaltyparb$ and $\next{\sduald}= \next{\duald}/\penaltypard$ leads to
\begin{subequations}
\begin{align}
\next{\m} &= \arg\min_{\m}  \reg(\m) +\frac{\penaltyparb}{2} \|\prev{\sdualb}-\A(\m)\prev{\u}\|_2^2 \label{msub}, \\
\next{\u} &= \arg\min_{\u} \frac{\penaltyparb}{2} \|\prev{\sdualb}-\A(\next{\m})\u\|_2^2 + \frac{\penaltypard}{2} \|\prev{\sduald} -\P\u\|_2^2,\label{usub}\\
\next{\sdualb}& = \prev{\sdualb} + \b -\A(\next{\m})\next{\u}, \label{dualb} \\
\next{\sduald}& = \prev{\sduald} + \d -\P\next{\u}, \label{duald}
\end{align}
\end{subequations}
beginning with $\sdualb_0=\b$ and $\sduald_0=\d$ and iterating until some stopping criterion is satisfied. In equations \eqref{msub} and \eqref{usub}, the original right-hand sides $\b$ and $\d$  of the wave equation and the observation equation are augmented with the running sum of the wave-equation and observation-equation errors as highlighted by the gradient ascent steps of Lagrange multipliers in equations~\eqref{dualb}-\eqref{duald}.


The $\m$-subproblem is a regularized least squares problem that can be solved efficiently for different forms of the regularization function \cite{Aghamiry_2019_IBC,Aghamiry_2019_CRO}.

Solving the $\u$-subproblem is however more challenging as it requires us to solve an augmented wave equation
\begin{equation} \label{AugWE}
\begin{pmatrix}
\sqrt{\penaltyparb}\A(\next{\m})\\
\sqrt{\penaltypard}\P
\end{pmatrix}
\next{\u}
=
\begin{pmatrix}
\sqrt{\penaltyparb}\prev{\sdualb}\\
\sqrt{\penaltypard}\prev{\sduald}
\end{pmatrix}
\end{equation}
in a least-squares sense. Although the normal operator $\penaltyparb \A(\next{\m})^T \A(\next{\m}) + \penaltypard \P^T \P$ associated with \eqref{AugWE}  is Hermitian positive definite, it has a wider bandwidth, a larger condition number and is denser than the indefinite Helmholtz operator $\A(\next{\m})$. These specifications can limit the applicability of direct and iterative solvers for large scale problems.\\  
Recently, \cite{Gholami_2022_EFW} showed that the exact solution of \eqref{usub} satisfies 
\begin{equation} \label{DAW}
 \next{\A} \next{\u}= \prev{\sdualb} + \next{\bold{S}}^{T}\next{\delta \d}^e,
\end{equation}
where $ \next{\A}\equiv \A(\next{\m})$, $\next{\bold{S}}=\P\next{\A}^{-1}$ is the modeling operator and
\begin{equation} \label{delta_de}
\next{\delta \d}^e = \left( \next{\bold{S}}\next{\bold{S}}^T  + \frac{\penaltyparb}{\penaltypard} \I\right )^{-1} \next{\delta \d}^r,
\end{equation}
with
\begin{equation}
\next{\delta \d}^r = \next{\d} - \next{\bold{S}} \next{\b}.
\end{equation}
The frequency domain implementation of the IR-WRI formulation proposed by \cite{Gholami_2022_EFW} is summarized in Algorithm \ref{Alg1}. 
\begin{algorithm}[htb!]
\caption{
New face of IR-WRI adapted for large scale problems.}
\label{Alg1}
\scriptsize
{\fontsize{8}{8}\selectfont
\begin{algorithmic}[1]
\STATE Begin with $k=0$ and an initial model $\bold{m}_0$.
\STATE Set ${\sdualb}_0=\bold{b}$ and ${\sduald}_0=\bold{d}$.
\item[]
\WHILE {convergence criteria  not satisfied}
\STATE $\next{\bold{S}}^T=\next{\A}^{-T}\P^T$
\item[]
\STATE $\next{\delta \d}^r = \next{\d} - \next{\bold{S}} \next{\b}$. 
\item[]
\hspace*{20.7em}%
        \rlap{\smash{$\left.\begin{array}{@{}c@{}}\\{}\\{}\\{}\\{}\\{}\\{}\\{}\end{array}\color{black}\right\}
          \color{black}\begin{tabular}{l}The new form of wavefield reconstruction instead of solving \eqref{AugWE}. \end{tabular}$}} 
\STATE $\next{\delta \d}^e = \left( \next{\bold{S}}\next{\bold{S}}^T  + \frac{\penaltyparb}{\penaltypard} \I\right )^{-1} \next{\delta \d}^r$
\item[]
\STATE $  \next{\u}= \next{\bold{A}}^{-1}[\prev{\sdualb} + \next{\bold{S}}^{T}\next{\delta \d^e}]$
\item[]
\STATE Update model parameters $\next{\m}$ by solving \eqref{msub}. 
\item[]
\STATE $\next{\sdualb} = \prev{\sdualb} + \b -\A(\next{\m})\next{\u}$
\item[]
\STATE $\next{\sduald} = \prev{\sduald} + \d -\P\next{\u}$
\item[]
\STATE $k = k+1$ ,
\item[]
\ENDWHILE
 \end{algorithmic}
}
\end{algorithm}
An issue associated with Algorithm \ref{Alg1} is the solution of the dense $\Nr \times \Nr$ normal system in the data space (step 6 and equation~\ref{delta_de}). This system can be easily solved in the frequency domain with direct solver suitable for dense linear systems once the normal operator $\next{\bold{S}}\next{\bold{S}}^T$ has been built explicitly. 
The adjoint of $\bold{S}$ can be built by performing $\Nr$ adjoint simulations in step 4 while the wavefield reconstruction in step 7 requires $\Ns$ simulations.
In summary, each iteration of the new frequency-domain IR-WRI algorithm requires $\Nr + \Ns$ wave-equation solutions to reconstruct the wavefields, while the classical algorithm requires $\Ns$ solutions of the augmented wave-equation in \eqref{AugWE}. These additional $\Nr$ simulations are a severe issue since one acquisition dimension is generally oversampled in stationary-recording experiments. This implies that $\Nr >> \Ns$ keeping in mind that receivers can be processed as sources and vice versa by virtue of the spatial reciprocity of the Green functions to mitigate the number of wave simulations during FWI. 

It is worth mentioning that the computational burden of the steps other than steps 4 and 7 in  Algorithm \ref{Alg1} is negligible. 
The next section reviews a sketching strategy to mitigate the computational burden of these multi-rhs simulations.

%
%


\subsection{Efficient multi-source/receiver processing with sketching methods} \label{sketch}
Random source encoding is a common strategy in imaging to decrease the computational burden of multi-source simulations. It simply consists in gathering several encoded sources into one or a few super-sources through a weighted summation. The drawback is that the FWI gradient is affected by an additional term representing cross-talk artefacts between each sources of the super-sources. However, this cross-talk term vanishes asymptotically when the random codes $\alpha_i$ satisfy
\begin{equation}
\mathbb{E}\left[\alpha_i^* \alpha_j\right] = \delta_{i,j}
\end{equation}
where $\mathbb{E}$ stands for the expectation over $\alpha$, $\delta_{i,j}$ is the Kronecker delta and the subscript of $\alpha$ labels one source of a super source \cite{Krebs_2009_FFW,Castellanos_2015_FFW1}. \\
\cite{Aghazade_2021_RSS} recast the concept of random source encoding in the more general framework of sketching methods and interface this approach with IR-WRI. We apply this sketching method to reduce the number of simulations during steps 4 and 7 of Algorithm \ref{Alg1}, which involve $\Nr$ and $\Ns$ right-hand sides, respectively. Accordingly, we design two sketching matrices $\bold{X} \in \mathbb{C}^{\Nr \times \Nr'}$ ($\Nr' \ll \Nr$) and $\bold{Y} \in \mathbb{C}^{\Ns\times \Ns'}$ ($\Ns' \ll \Ns$), which satisfy
\begin{equation}
\mathbb{E}\left[\bold{X} \bold{X}^T\right] = \bold{I},   ~~~ \mathbb{E}\left[\bold{Y} \bold{Y}^T\right] = \bold{I}.
\end{equation}
Among the possible choices reviewed in \cite{Aghazade_2021_RSS}, we use Gaussian sketching matrices, which means that the elements of $\bold{X}$ and $\bold{Y}$ are independently sampled from the standard normal distribution.
Then, we right multiply the right hand side of the equations involved in step 4 and 7 by $\bold{X}$ and $\bold{Y}$, respectively, to decrease the number of PDE solutions from $\Nr$ to $\Nr'$ and from $\Ns$ to $\Ns'$, respectively. Moreover, we regenerate the projection matrices $\bold{X}$ and $\bold{Y}$ at each IR-WRI iteration to remove more efficiently the cross-talk noise.

\subsection{The convergent Born series method} \label{CBSsec}

\subsubsection{Principles}
We solve the Helmholtz equation efficiently with the iterative convergent Born series (CBS) method \citep{Osnabrugge_2016_CBS}, whose main principles are reviewed here. The reader is referred to \cite{Osnabrugge_2016_CBS} for more details. \\
Let us decompose $\m$ as the sum of a homogeneous background $\c_0$ and perturbation model $\delta\m$, $\m=\c_0+\delta\m$, then \eqref{eq:PDE} can be written as
\begin{equation} \label{eq:PDE2}
\del \u+ \w^2\Diag(\c_0)\u=\b - \w^2\Diag(\delta\m)\u,
\end{equation}
which may be solved with the following fixed-point iteration
\begin{equation} \label{eq:PDE_fp}
\nextin{\u} = \A_0^{-1}(\b - \w^2\Diag(\delta\m)\previn{\u}),
\end{equation}
where 
\begin{equation}
\A_0^{-1} =  \left(\del + \w^2\Diag(\c_0)\right)^{-1}.
\end{equation}
Beginning with $\u_0=\A_0^{-1}\b$, the first iteration of this series is known as the first order Born series.
The series \eqref{eq:PDE_fp} is convergent as long as the spectral radius of $\w^2\A_0^{-1}\Diag(\delta\m)$ is less than unity which is the case for small perturbations $\delta\m$. Otherwise, the series will diverge and we need an appropriate pre-conditioner to make it convergent.

Following \cite{Osnabrugge_2016_CBS}, we first use the Green's functions corresponding to an attenuative medium by modifying equation \eqref{eq:PDE2} as
\begin{equation} \label{eq:PDE2_atten}
\del \u+ \Diag(\w^2\c_0-\iu \epsilon)\u=\b  - \Diag(\w^2\delta\m+\iu \epsilon)\u,
\end{equation}
with the imaginary unit $\iu$, leading to the fixed-point problem
\begin{equation} \label{eq:PDE_fp_atten}
{\u} = \aA^{-1}(\b - \Diag(\w^2\delta\m+\iu\epsilon){\u}),
\end{equation}
where
\begin{equation}
\aA^{-1} =  \left(\del + \Diag(\w^2\c_0-\iu \epsilon)\right)^{-1},
\end{equation}
denotes the attenuated Green's functions operator. Then multiplying a properly defined pre-conditioner
\begin{equation}
\M = \frac{-\iu}{\epsilon}\Diag(\w^2\delta\m+\iu\epsilon)
\end{equation}
to both sides of \eqref{eq:PDE_fp_atten}, after adding $\u$ to both sides and rearranging the terms, gives the following modified Born series
\begin{equation} \label{CBS}
\nextin{\u}= \previn{\u} - \M \previn{\u} + \M \aA^{-1}(\b - \Diag(\w^2\delta\m+\iu\epsilon) \previn{\u}).
\end{equation}
\cite{Osnabrugge_2016_CBS} showed that this series is convergent for
\begin{equation} \label{epsi}
\epsilon \geq \w^2\max(|\delta\m|).
\end{equation}

The Green's function of the attenuative homogeneous medium can easily be found in the Fourier domain as
\begin{equation} \label{Greens_anal}
\aA^{-1} =  \left(\del + \Diag\left(\w^2\c_0-\iu\epsilon\right)\right)^{-1} = \F^T \Diag\left(\frac{1}{-|\bold{k}|^2 + \w^2\c_0-\iu\epsilon}\right)\F,
\end{equation}
where $\F$ is the 2D spacial Fourier transform and $\bold{k}$ denotes the wavenumber vector. 

The CBS iteration in \eqref{CBS} can be simplified as
\begin{align} \label{CBS2}
\nextin{\u} &= \previn{\u} - \M \aA^{-1}(\A\previn{\u}-\b)\\
&= \previn{\u} - \bold{H}^{-1}\A^T(\A\previn{\u}-\b),
\end{align}
where $\bold{H}=\A^T\aA\M^{-1}$. Thus this iteration is basically a Newton-type method applied for minimizing $\|\A\u-\b\|_2^2$ by an approximate Hessian $\bold{H}=\A^T\aA\M^{-1}\approx \A^T\A$. However, this is different from the recent backward-forward time-stepping recursion proposed to directly solve the normal system for the wavefields \citep{Aghamiry_2019_AEW}.

\subsubsection{Numerical complexity of the CBS method}
\cite{Osnabrugge_2016_CBS} show that the CBS method expands wave packet at a pseudo speed $\Delta r$ given by
\begin{equation}
\Delta r = \frac{2 k_0}{\epsilon}.
\label{eqdr}
\end{equation}

Accordingly, the convergence speed is optimized by minimizing $\epsilon$, which gives for $k_0$
\begin{equation}
k_0^2=\left( \left(\text{Re}\left\{k(\bold{x})^2\right\}\right)_m+\left(\text{Re}\left\{k(\bold{x})^2\right\}\right)_M \right)/2,
\label{eqk0}
\end{equation}
where $(\cdot)_m$ and $(\cdot)_M$ are the minimum and maximum values of vector $(\cdot)$ and $\text{Re}\left\{\cdot\right\}$ is the real part of the complex number $(\cdot)$.

Considering that the time complexity of one iteration of the CBS method scales to $\mathcal{O}(n \log n) \equiv \mathcal{O}(k_M^D \log k_M^D)$ where $D \in \{1,2,3\}$ is the dimensionality of the simulation, the time complexity of the full run is obtained by dividing the above complexity by the pseudo speed (number of propagated wavelength per iteration) leading to $\mathcal{O}(\nu  k_M^{D+1} \log k_M^D)$ where $\nu = \frac{k_M}{k_0}-\frac{k_0}{k_M}$ is the scattering contrast \citep[][ Table 1]{Osnabrugge_2016_CBS}. The expression of the complexity shows that the convergence speed of the CBS method is controlled both by frequency and scattering contrast. Note that the above complexity is provided for a grid interval set to half the minimum wavelength. 

\section{Numerical results}
\begin{figure}[htb!]
\center
\includegraphics[width=1\columnwidth,trim={0 0cm 0 0cm},clip]{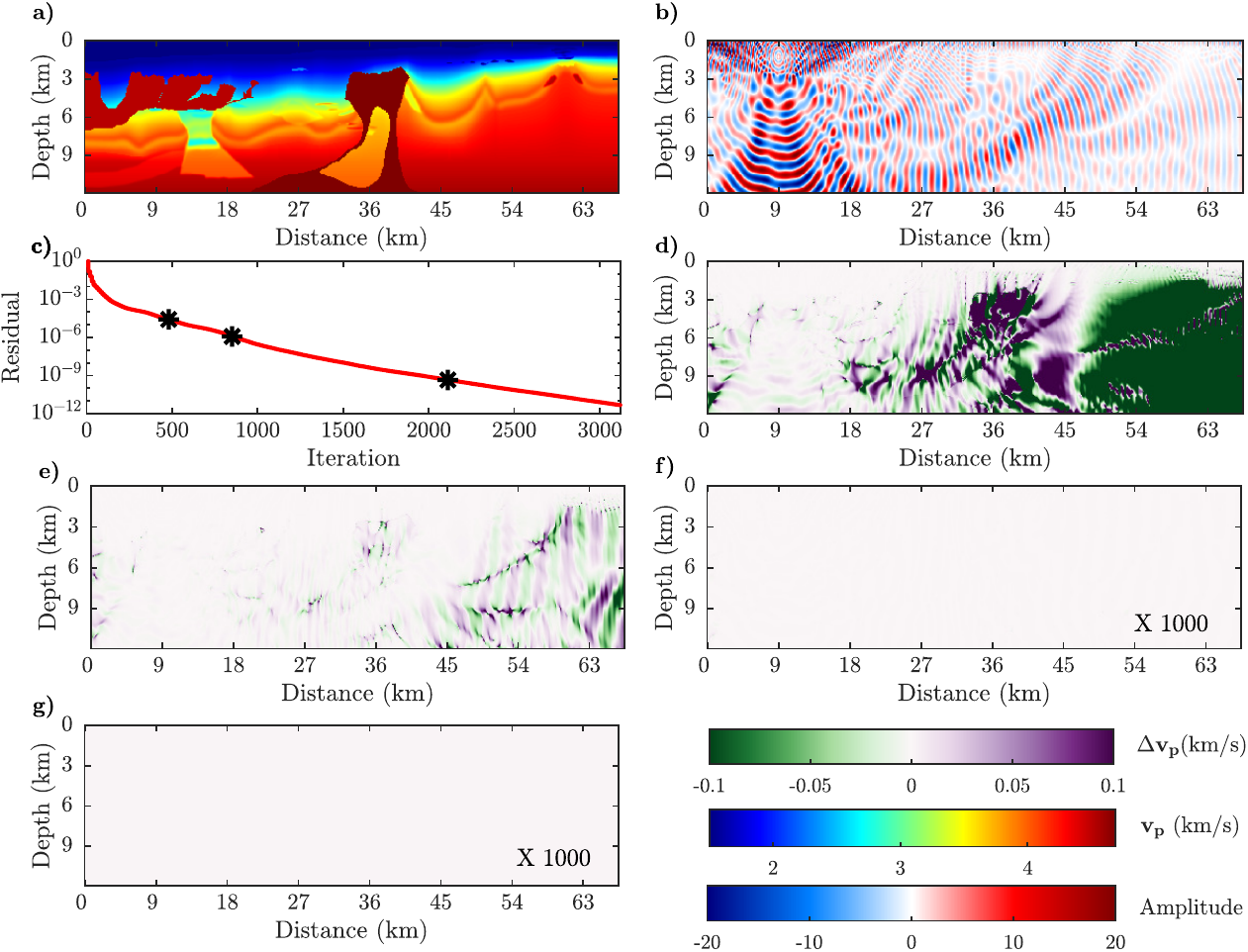}
\caption{(a) 2004 BP salt model. (b) Computed 3Hz wavefield with CBS for $\eta$=1e-8. 
A linear gain with distance is applied from the source for amplitude balancing. (c) Relative wave-equation error as a function of iteration count for $\eta$=1e-10. The three stars point the number of iterations required to reach $\eta$=1e-3, $\eta$=1e-5, and $\eta$=1e-8. (d) Reconstructed velocity model errors from the 3 Hz wavefields computed with $\eta$ = (d) 1e-3, (e) 1e-5, (f) 1e-8, and (g) 1e-10. Velocity models are reconstructed from simulated wavefieds according to equation \ref{eq:PDEsf}.}
\label{Fig_BP_forward}
\end{figure}

We first assess the accuracy of the CBS forward engine with the 2004 BP salt model (Figure \ref{Fig_BP_forward}a) \cite{Billette_2004_BPB}. 
Then, we benchmark the CBS-based IR-WRI implemented without and with source sketching(Algorithm~\ref{Alg1}) against the original formulation of IR-WRI where the augmented wave-equation is solved with a finite-difference (FD) method and a direct solver. 
We set $k_0$ according to eq.~\ref{eqk0} to perform CBS simulations.
We implement the original formulation of IR-WRI with a 9-point FD staggered-grid stencil with PML boundary condition and anti-lumped mass to solve the Helmholtz equation \cite{Chen_2013_OFD}. Moreover, we set $\lambda$ as a small fraction of the largest eigenvalue of $\bold{A}_0^{-T}\bold{P}^T\bold{P}\bold{A}_0^{-1}$ to perform IR-WRI \cite{vanLeeuwen_2016_PMP}. 
\begin{figure}[htb!]
\center
\includegraphics[width=1\columnwidth,trim={0 0cm 0 0cm},clip]{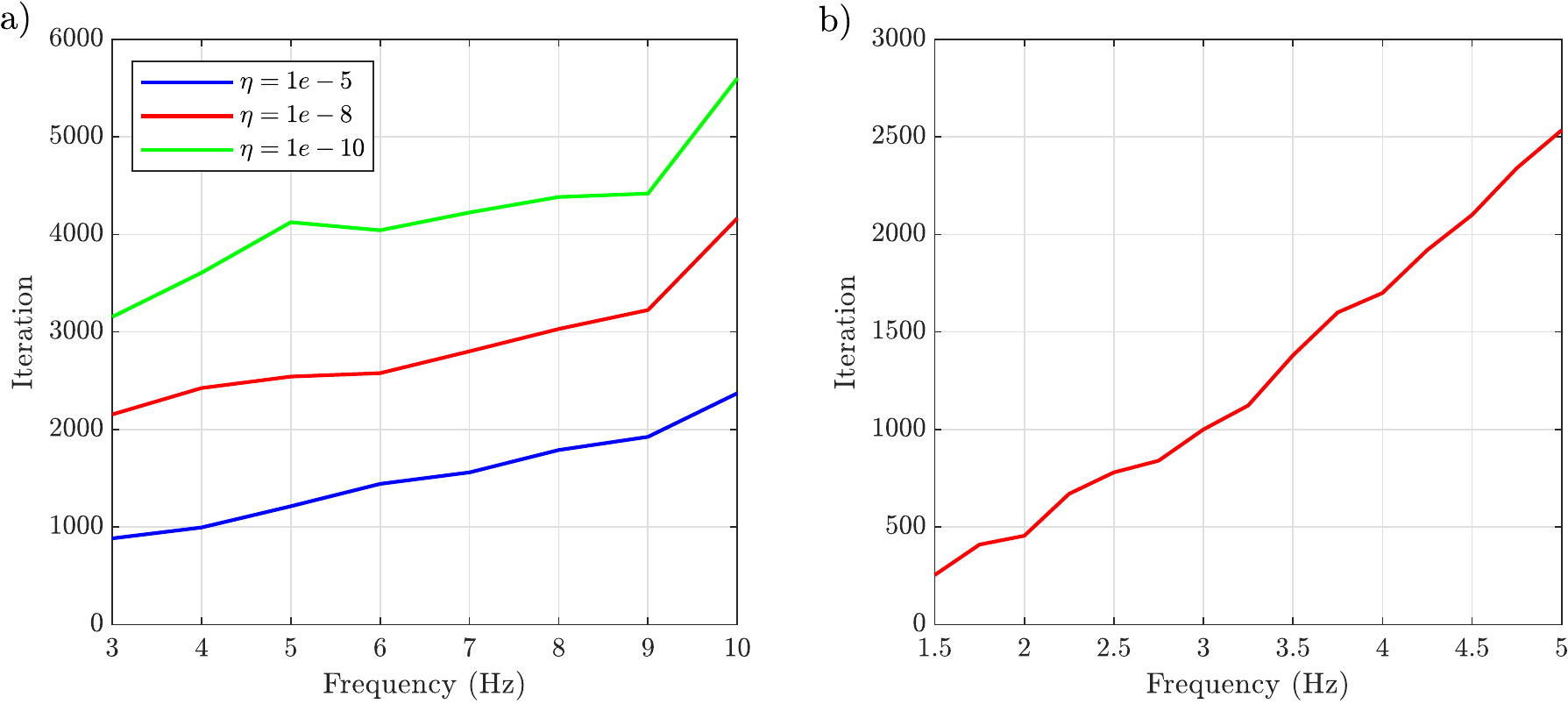}
\caption{Number of iterations required to compute wavefields with CBS in 2004 BP salt model as a function of frequencies. (a) Iteration counts with true velocity model for $\eta=$1e-5 (blue), $\eta=$1e-8 (red), and $\eta=$1e-10 (green). (b) Iteration counts with the final IR-WRI velocity models of each mono-frequency inversion extracted in section \ref{invres} for $\eta=$1e-8.}
\label{Fig_BP_scalability}
\end{figure}
\subsection{Forward modeling with CBS and stopping criterion of iteration}
We define the stopping criterion of iteration that guarantees accurate wavefield reconstruction with the CBS method as
\begin{equation}
\delta \b=\frac{\| \bold{A} \u - \b\|_2}{\|\b\|_2} \le \eta,
\end{equation}
where $\delta \b$ is the relative wave-equation error and $\eta$ is a user-defined threshold.
To assess the accuracy of the wavefield as a function of $\eta$, we use the bilinearity of the wave-equation that allows one to reconstruct exactly the model parameters from a known monochromatic wavefield by a pointwise division as
\begin{equation} \label{eq:PDEsf}
\m=\frac{\b-\del \u}{\w^2\u},
\end{equation}
where the above equality can be easily inferred from Eq. \ref{eq:PDE} and the Laplacian in the right-hand side is computed in the Fourier domain to prevent discretization error. Any inaccuracies in the simulated wavefields will translate into artefacts in the computed velocity model with equation ~\ref{eq:PDEsf}, hence providing an indirect assessment of the wavefield accuracy. \\
We compute a 3-Hz monochromatic wavefield using a point source at (1.5, 9) \textit{km} and absorbing boundary conditions along all the model's edges. Figure \ref{Fig_BP_forward}b shows the reconstructed wavefield for $\eta=$1e-10. Figures \ref{Fig_BP_forward}(d-g) show the error in the reconstructed velocity models with eq.~\ref{eq:PDEsf} when the velocity models are built from the computed wavefields with $\eta$=1e-3, 1e-5, 1e-8 and 1e-10, respectively. One can readily check that $\eta$=1e-8 and $\eta$=1e-10
lead to an accurate reconstruction of the BP salt model, while the arterfacts shown in the velocity models inferred from the computed wavefields with $\eta$=1e-3 and $\eta$=1e-5, Figures \ref{Fig_BP_forward}(d-e), highlight the insufficient accuracy of these wavefields.  Hereafter, we use $\eta=1e-8$ to compute wavefields with CBS except where expressly stated otherwise.
Also, the relative wave-equation error $\delta \bold{b}$ for this test is plotted as a function of iteration number in Figure \ref{Fig_BP_forward}c.\\ 
\subsection{Numerical convergence of CBS method}
We continue by assessing the two factors that control the rate of the convergence of CBS, namely the frequency and the scattering potential (see eq.~\ref{eqdr}). 
Eq. \eqref{eqdr} shows that the number of CBS iterations increases with frequencies for a fixed subsurface model. We verify numerically this statement for different values of $\eta$ and the true 2004 BP salt model in Figure \ref{Fig_BP_scalability}a, which shows that the number of iterations required to reach a specific $\eta$ grows approximately linearly with frequency.
Moreover, according to Eq. \eqref{eqdr}, models with large contrasts require more CBS iterations for convergence than smoother models \cite{Osnabrugge_2016_CBS}.  
When considering a multiscale waveform inversion with frequency continuation, the contrasts between the reconstructed models by IR-WRI and the homogeneous background model used by CBS become sharper as the frequency involved in the inversion increases. Accordingly, the number of iterations of CBS during the early stages of the inversion (i.e., at low frequencies and for smooth models) is smaller than that required during the final stages (i.e., at higher frequencies and for sharper models). This is illustrated in Figure \ref{Fig_BP_scalability}b with the reconstructed velocity models by IR-WRI presented in the next section. Figure \ref{Fig_BP_scalability}b shows how the number of CBS iteration increases with frequency when the simulations are performed in the final IR-WRI velocity models of each mono-frequency inversion with the frequency involved in the inversion. Accordingly, the figure highlights the computational efficiency of the CBS method at low frequencies where a limited number of iterations are necessary to reach the convergence criterion.
\begin{figure}[htb!]
\center
\includegraphics[width=1\columnwidth,trim={0 0cm 0 0cm},clip]{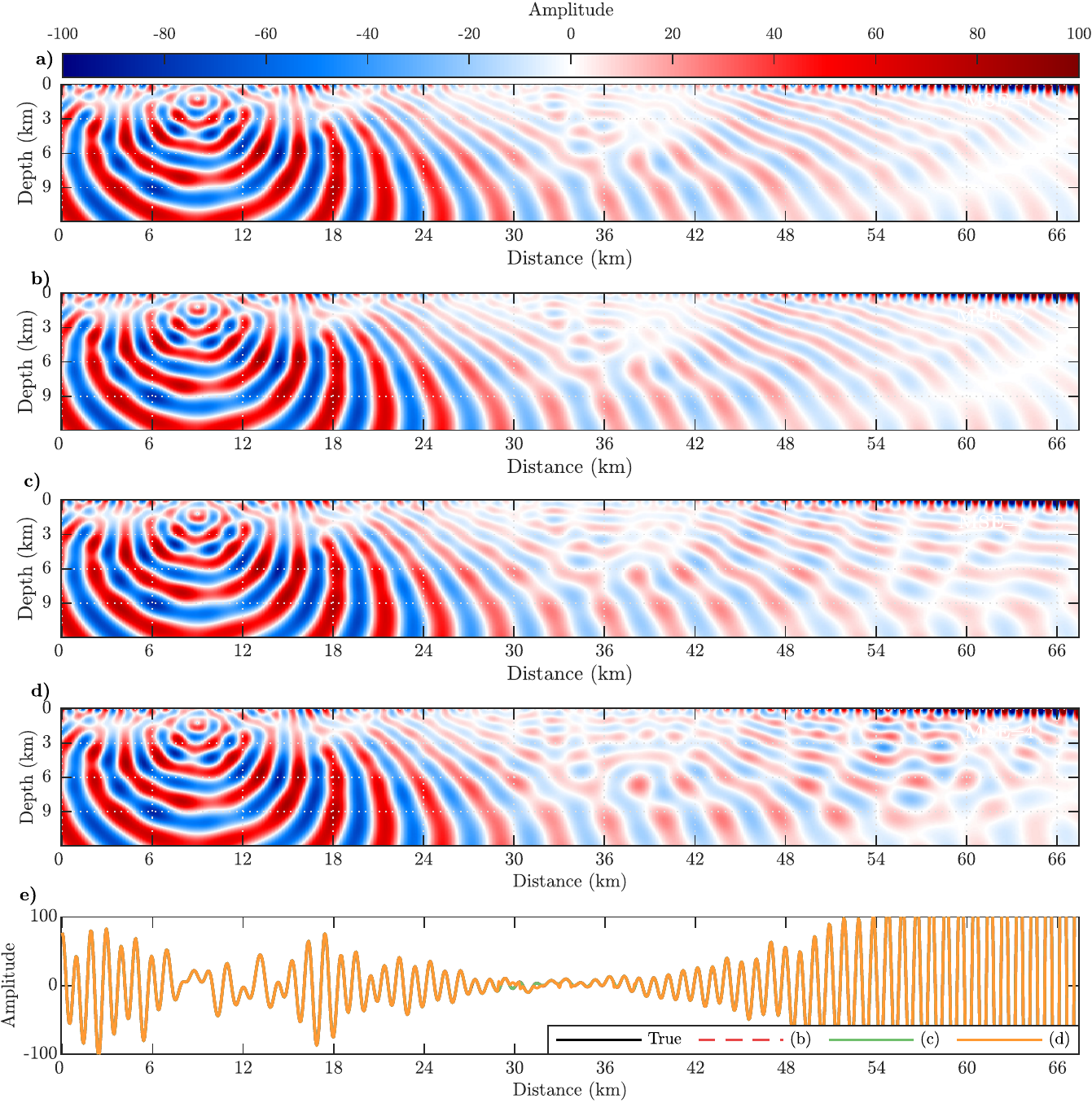}
\caption{Reconstructed DA wavefields (a) using the FD method and a direct solver (eq. \ref{eq:PDEsf}), (b-d) with CBS (eq. \ref{DAW})  (b) without receiver sketching and (c-d) with receiver sketching when the number of blended receivers is set to (c) 200 (77\% reduction) and (d) 100 (88\% reduction). (e) The true wavefield and the reconstructed ones shown in (a-d) are superimposed at receiver positions. A linear gain function with distance from the source is applied on all the panels for amplitude balancing.}
\label{DA_wavefield15}
\end{figure}
\subsection{Computing DA wavefield with CBS and with/without receiver sketching}
We now assess the CBS method to reconstruct DA wavefields with eq. \ref{DAW}. We use a point source at (1.5, 9) \textit{km} and a line of 900 receivers spaced 75~m apart at the surface. 
We compute the  DA wavefields in a crude laterally-homogeneous velocity-gradient model with velocities ranging between 1.5 to 3.5 km/s while the recorded data are computed in the true BP salt model. Note that we use $\eta$=1e-10 to compute the recorded data, while we use $\eta$=1e-8 to compute the DA wavefield, which means that we don't use an inverse crime procedure to perform DA wavefield reconstruction  with CBS. As benchmark solution, we compute the 1.5 Hz DA wavefield by solving the normal equation, eq. \ref{AugWE}, with a sparse direct solver when the wave-equation operator is discretized with the 9-point FD method (Figure \ref{DA_wavefield15}a).
The computed DA wavefield with CBS (from eq. \ref{DAW}) and without receiver sketching is shown in Figure \ref{DA_wavefield15}b and closely match the FD-based DA wavefield shown in Figure \ref{DA_wavefield15}a, which provides a first validation of the implementation of eq.~\ref{DAW}. \\
One difficulty with eq.~\ref{DAW} is related to the computation and storage of $\next{\bold{S}}^T$ (line 4 of Algorithm \ref{Alg1}). Receiver sketching was proposed in section \ref{sketch} as a remedy for decreasing the dimension of $\bold{S}$ and hence the number of receiver-side simulations. 
We now assess how receiver sketching impacts the accuracy of DA wavefields compared to the case when sketching is not used. DA wavefields reconstructed with 200 sketched (or super) receivers and 100 sketched receivers are shown in Figure \ref{DA_wavefield15}(c-d), leading to a reduction of the number of receivers by a factor $\sim$4 and 8, respectively. At short offsets, the extracted wavefields without and with receiver sketching are similar, while, we show a continuous degradation of the accuracy of DA wavefields as a function of distance from the source.
The footprint of this accuracy degradation on the convergence of IR-WRI is assessed in the next section. 
Finally, the direct comparison between the true wavefields and the DA wavefields at receiver positions show that the match is almost perfect for all the cases (Figures \ref{DA_wavefield15}(b-d)). But, this statement is not verified at higher frequencies, an issue that will be discussed at the end of this section. \\
\begin{figure}[htb!]
\center
\includegraphics[width=1\columnwidth,trim={0 0cm 0 0cm},clip]{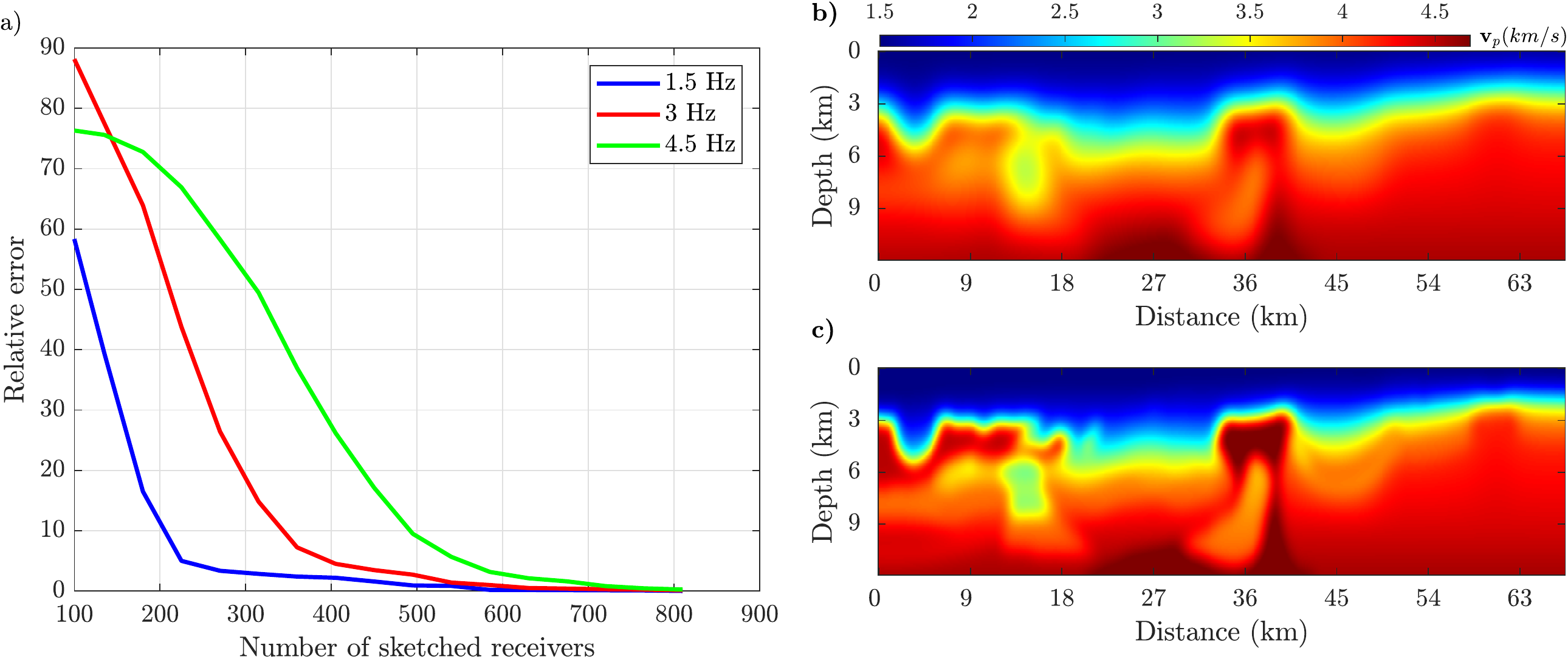}
\caption{(a) Relative error between computed DA wavefields with CBS, eq.~\ref{DAW}, without receiver sketching and with receiver sketching for 1.5 Hz (blue), 3 Hz (red), and 4.5 Hz (green). The horizontal coordinate shows the number of sketched receivers. For 1.5~Hz frequency, the laterally-homogeneous velocity-gradient model with velocities ranging between 1.5 to 3.5 km/s is used, while the velocity models (b) and (c) are used for 3 and 4.5 Hz simulations, respectively.}
\label{Sketching_error}
\end{figure}
\subsubsection{Impact of receiver sketching with frequency}
\begin{figure}[htb!]
\center
\includegraphics[width=1\columnwidth,trim={0 0cm 0 0cm},clip]{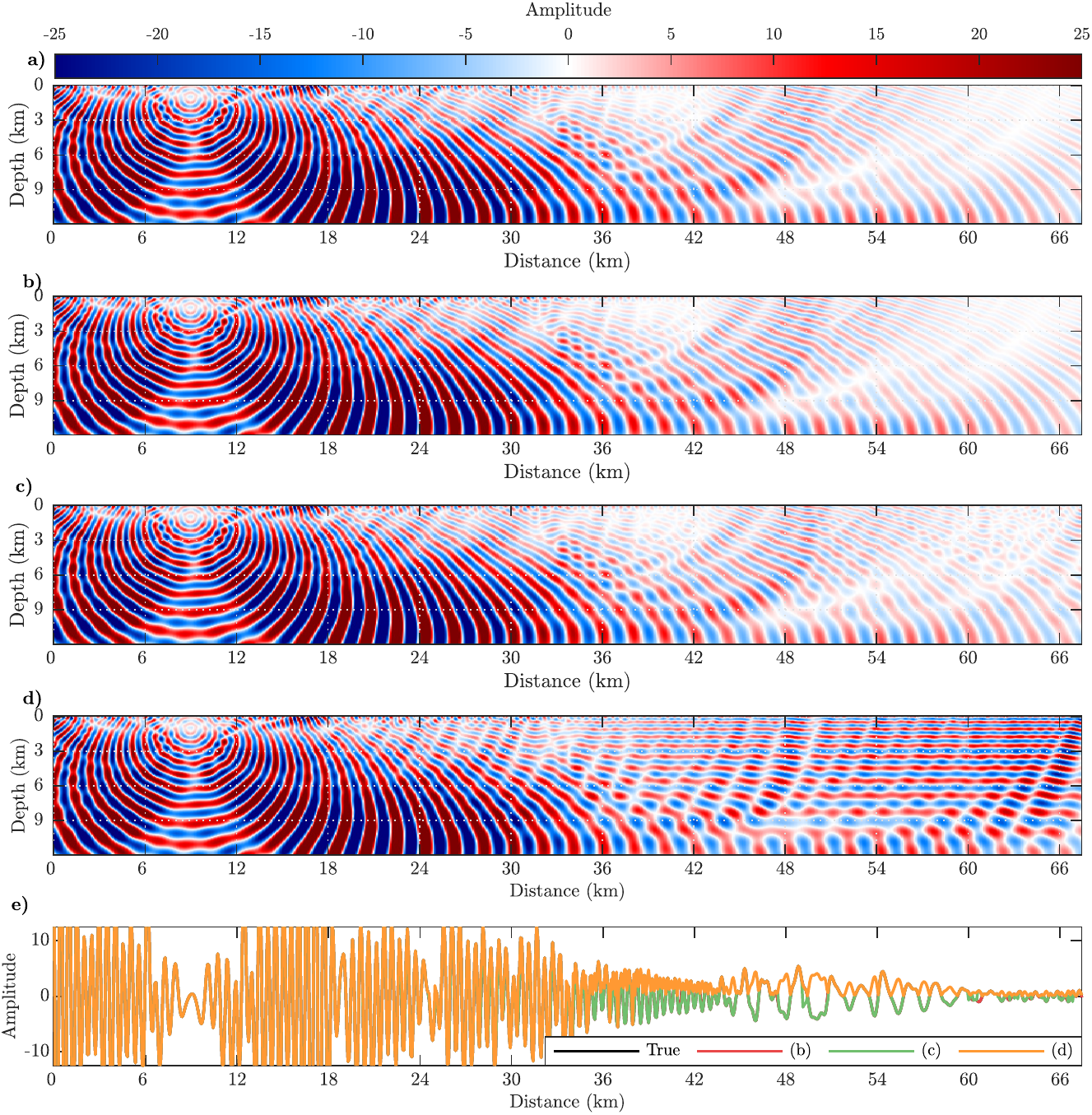}
\caption{Same as Figure \ref{DA_wavefield15}, but for 3 Hz and the initial model depicted in Figure \ref{Sketching_error}b. Also panel (c) is with 300 sketched receivers (66\% reduction), and (d) is with 200 (77\% reduction). }
\label{DA_wavefield30}
\end{figure}

\begin{figure}[htb!]
\center
\includegraphics[width=1\columnwidth,trim={0 0cm 0 0cm},clip]{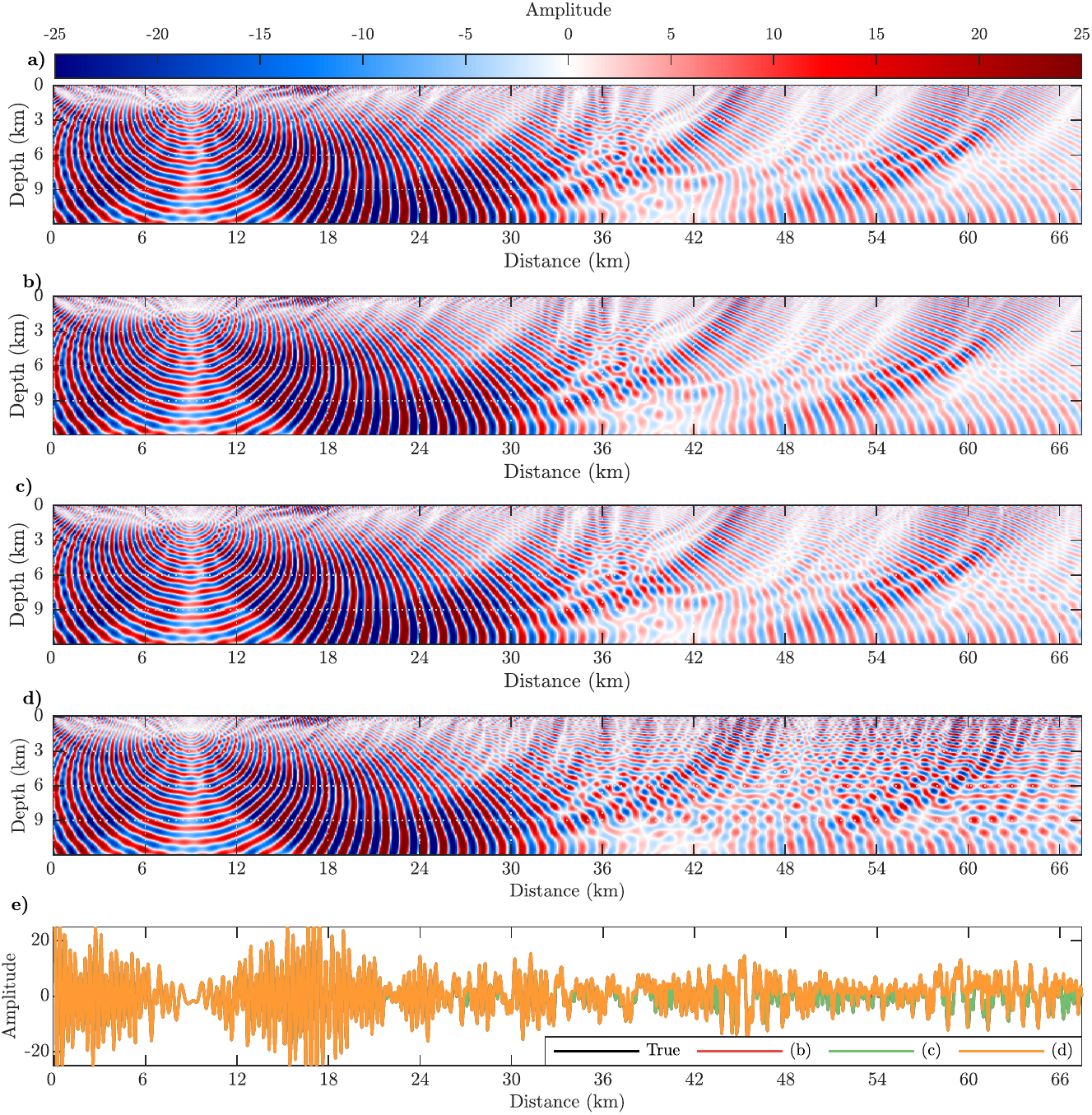}
\caption{Same as Figure \ref{DA_wavefield15}, but for 4.5 Hz and the initial model depicted in Figure \ref{Sketching_error}c. Also panel (c) is with 400 sketched receivers (55\% reduction), and (d) is with 300 (66\% reduction).}
\label{DA_wavefield45}
\end{figure}

The footprint of receiver sketching depends not only on the number of blended receivers but also on the frequency and closeness of the initial model to the true model. To investigate further this footprint, we plot in Figure \ref{Sketching_error}a the relative error between the reconstructed DA wavefield with different numbers of sketched (or super) receivers and the reference wavefield computed without sketching (Figure \ref{DA_wavefield15}b) for the 1.5 Hz frequency in blue  when the initial velocity gradient model is used. We also show the same curve for the 3 Hz and 4.5 Hz frequencies when the models depicted in Figures \ref{Sketching_error}b-\ref{Sketching_error}c are used as initial models, respectively. These initial models are closer to the true model to mimic those that are reconstructed by FWI when increasing frequency. We show that the wavefield error increases as the number of sketched receivers decreases but at a different rate as a function of frequency.
To have a better visual understanding, we repeat the Figure \ref{DA_wavefield15}, but for 3 Hz with the initial model depicted in Figure \ref{Sketching_error}b, and for 4.5 Hz with the initial model depicted in Figure \ref{Sketching_error}c. We need more sketched receivers to build wavefield with acceptable accuracy when the frequency increases. Accordingly, we increase the number of sketched receivers in Figures \ref{DA_wavefield30}c-\ref{DA_wavefield30}d to 300 and 200 and in Figures \ref{DA_wavefield45}c-\ref{DA_wavefield45}d to 400 and 300. Also, like Figures \ref{DA_wavefield15}(c-d), the continuous degradation of the accuracy of DA wavefields with sketched receivers as a function of distance from the source are obvious, especially for those that have a smaller number of sketched receivers (Figures \ref{DA_wavefield30}d and \ref{DA_wavefield45}d). Finally, the direct comparison between the true wavefields and the DA wavefields at receiver positions in these figures (Figures \ref{DA_wavefield30}e and Figures \ref{DA_wavefield45}e) show obviously the degradation of the DA wavefield accuracy at the long offsets as the number of sketched receivers decreases. 

\subsection{Inversion results} \label{invres} 
\begin{figure}
\center
\includegraphics[width=1\columnwidth,trim={0 0cm 0 0cm},clip]{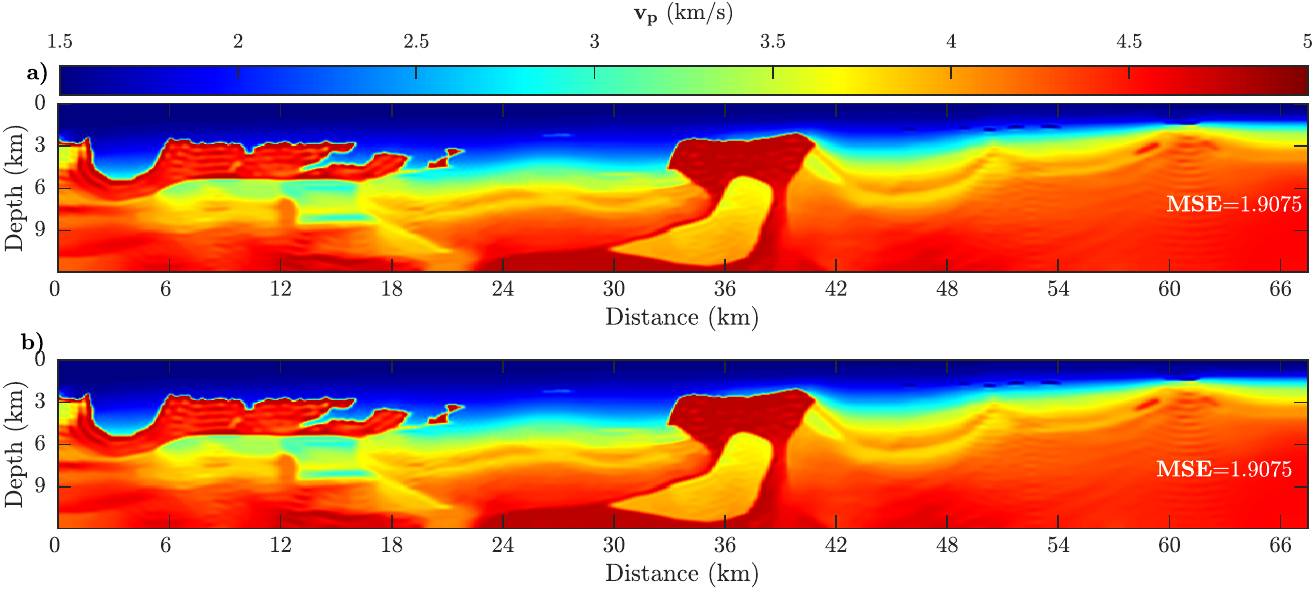}
\caption{IR-WRI results obtained with (a) the original formulation, (b) the new face formulation when the FD method is used as forward engine.}
\label{BP_inv_results_FD}
\end{figure}

\begin{figure}
\center
\includegraphics[width=1\columnwidth,trim={0 0cm 0 0cm},clip]{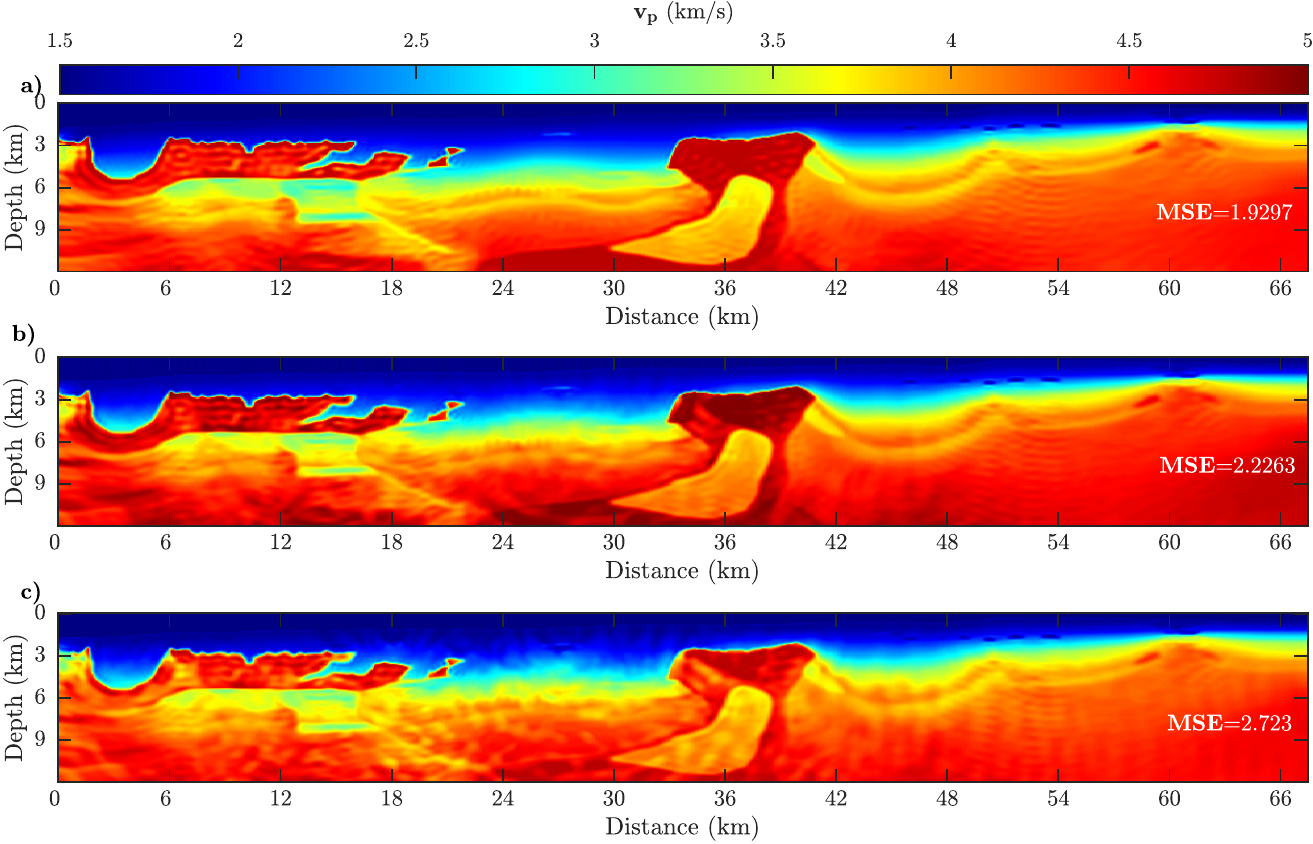}
\caption{New face IR-WRI results using CBS (a) without sketching, (b-c) with source/receiver sketching. The number of sketched receivers increases from 100 to 450 as frequencies increases from 1.5 Hz to 5 Hz. The number of sketched sources are fixed and are equal to 30 in (b) and 10 in (c).}
\label{BP_inv_results}
\end{figure}

\begin{figure}
\center
\includegraphics[width=1\columnwidth,trim={0 0cm 0 0cm},clip]{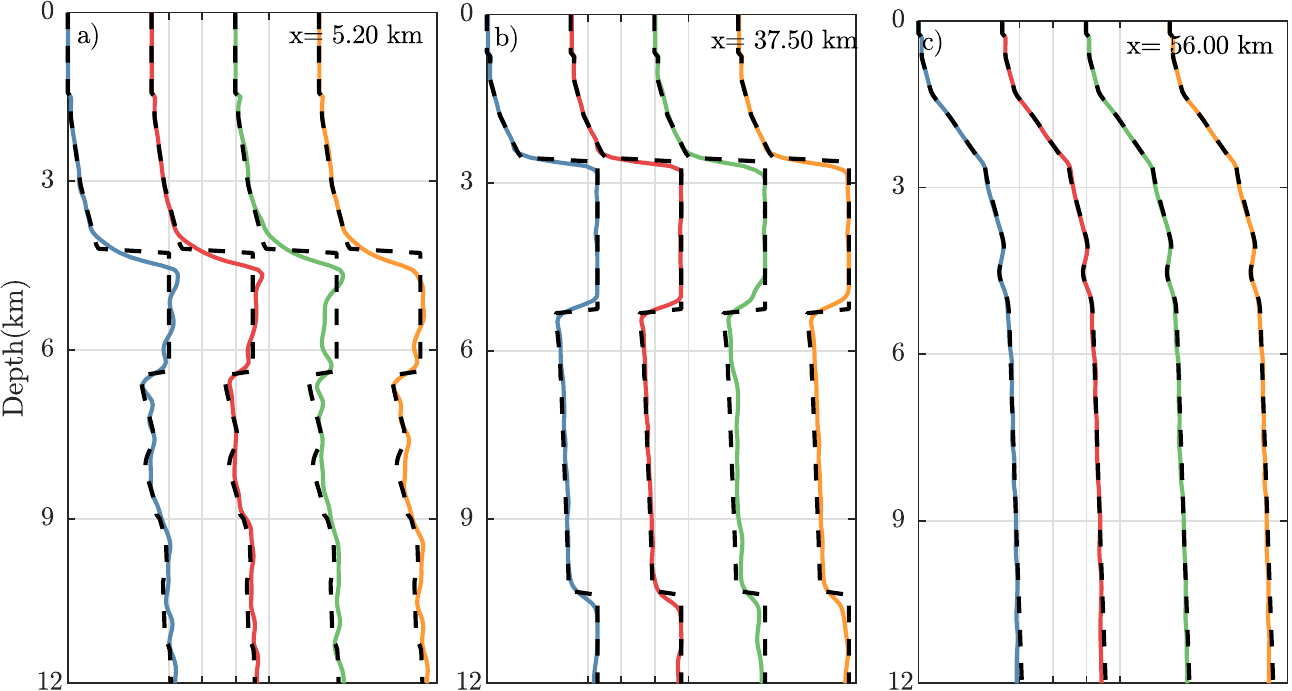}
\caption{A direct comparison between true (dashed black) and the reconstructed velocity of Figures \ref{BP_inv_results_FD}a (or \ref{BP_inv_results_FD}b) (blue), \ref{BP_inv_results}a (red), \ref{BP_inv_results}b (green), and \ref{BP_inv_results}c (orange)  at (a) $X$=5.2~km, (b) $X$=37.5~km, and (c) $X$=56.0~km.}
\label{BP_inv_results_logs}
\end{figure}

We now assess the new implementation of IR-WRI when the 1D gradient velocity model (laterally-homogeneous with velocities ranging between 1.5 to 3.5 km/s) is used as the starting model. 
The acquisition consists of 67 sources spaced 1~km apart on the seabed and 900 receivers spaced 75~m apart at the surface. The source signature is a 4~Hz Ricker wavelet. We apply the inversion on a long-offset fixed-spread dataset in the 1.5~Hz - 5~Hz frequency band with a frequency interval of 0.25~Hz. 
Following a multiscale frequency continuation strategy, we proceed over small frequency batches from low frequencies to higher ones when each batch of inversion contains two frequencies with one frequency overlap. We perform three paths through the frequencies, using the final model of one path as the initial model of the next one. The starting and finishing frequencies of the three paths are [1.5, 2], [1.5, 4], [3, 5]~Hz, respectively. The stopping criterion of iterations is a maximum of ten iterations per batch in all the cases. For all the inversion tests , we apply adaptive regularization and bound constraints \citep{Aghamiry_2020_FWI}. \\ 
We first check the equivalency between the inversion results of the original formulation of IR-WRI implemented with eq. \ref{AugWE} and the new face when the FD method is used as forward engine (Figure \ref{BP_inv_results_FD}). No receiver sketching is used at this stage. The results validate that both of the formulations reach the same minimizers when they are performed with the same experimental setup. \\
Then, we perform the new face IR-WRI using the CBS method as forward engine. The velocity model obtained without source and receiver sketchings (Figure \ref{BP_inv_results}a) is close to the models obtained with the FD forward engine (Figure \ref{BP_inv_results_FD}) albeit they are not identical. These differences probably results from the fact that we didn't use an inverse crime since the recorded data and the simulated ones are computed with two different stopping criteria of iterations ($\eta$=1e-10 and $\eta$=1e-8, respectively). 
Finally, we repeat the test with receiver and source sketching. According to the conclusions revealed by Figure \ref{Sketching_error}, we match the number of sketched receivers to frequency. We start with 100 sketched receivers at 1.5 Hz frequency, and finally use 450 sketched receivers at 5 Hz. On the source side, we keep the number of sketched sources fixed to 30 and 10 where the reconstructed models are shown in Figures \ref{BP_inv_results}b-\ref{BP_inv_results}c, respectively. We also show direct comparisons between the true velocity model and the reconstructed models shown in Figures \ref{BP_inv_results_FD} and \ref{BP_inv_results} at different distances in Figure \ref{BP_inv_results_logs}. 
The numerical results and the model errors, which are written in each panel of the figures, show that the extracted velocity model without source/receiver sketching (Figure \ref{BP_inv_results}a) and with 30 sketched sources (Figure \ref{BP_inv_results}b) are close together when the latter ones are extracted with 71 \% reduction in the number of PDE resolution and memory requirement. The inversion result with 10 sketched sources (Figure \ref{BP_inv_results}c) has a 73 \% reduction in the number of PDE resolution and memory requirements, but the inversion result has some high-frequency cross-talk artifacts because of the low number of sketched sources.
The details about the number of PDE resolutions are outlined in Table \ref{Tab1}.
\begin{table}[]
\begin{center}
\caption{The number of PDE resolutions for The reconstructed models using the new-face of IR-WRI represented in Figures \ref{BP_inv_results}b-\ref{BP_inv_results}d.}
\begin{tabular}{|c|c|c|c|l}
\cline{1-4}
\textbf{\begin{tabular}[c]{@{}c@{}}Number of PDE\\  resolution\end{tabular}} & \textbf{\begin{tabular}[c]{@{}c@{}}Without \\ sketching\end{tabular}} & \textbf{\begin{tabular}[c]{@{}c@{}}With receiver sketching \\ and 30 sketched sources\end{tabular}} & \textbf{\begin{tabular}[c]{@{}c@{}}With receiver sketching \\ and 10 sketched sources\end{tabular}} &  \\ \cline{1-4}
\textbf{Backward}                                                            & 180000                                                                & 51000 (72 \% reduction)                                                                            & 51000 (72 \% reduction)                                                                            &  \\ \cline{1-4}
\textbf{Forward}                                                             & 13400                                                                 & 6000 (66 \% reduction)                                                                             & 2000 (85\% reduction)                                                                              &  \\ \cline{1-4}
\textbf{Total}                                                               & 193400                                                                & 57000 (71 \% reduction)                                                                            & 53000 (73 \% reduction)                                                                            &  \\ \cline{1-4}
\end{tabular}
\label{Tab1}
\end{center}
\end{table}
%

\begin{figure}
\center
\includegraphics[width=1\columnwidth,trim={0 0cm 0 0cm},clip]{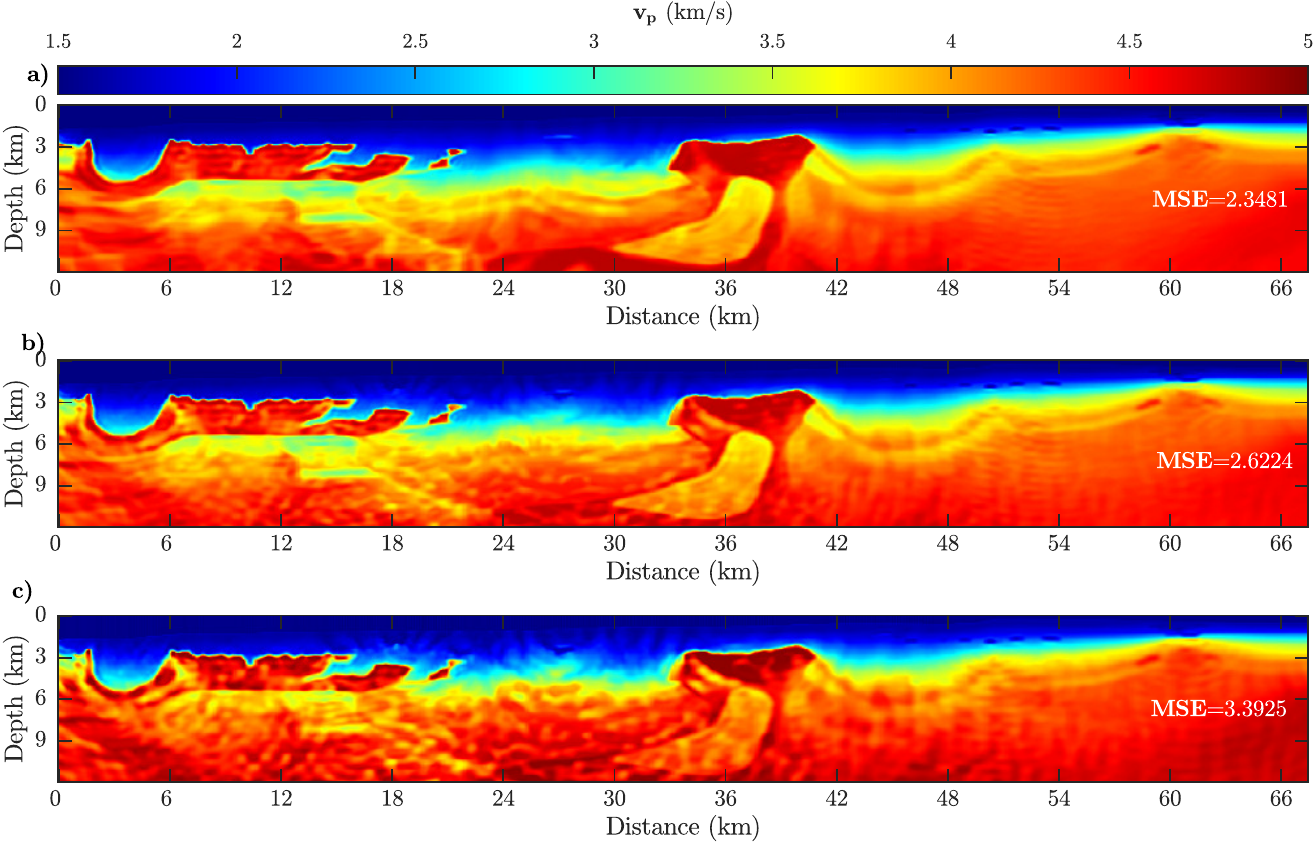}
\caption{New face IR-WRI results obtained with the FD forward engine. (a) No sketching is used. (b-c) Source/receiver sketching using the same configuration as in Figures \ref{BP_inv_results}b-\ref{BP_inv_results}d. The recorded data are computed with CBS to avoid inverse crime.}
\label{BP_inv_results_IC}
\end{figure}

One interesting feature of CBS is its high accuracy for solving the PDE as it is not affected by discretization errors \cite{Osnabrugge_2016_CBS}. As the final test, we generate a dataset using CBS with $\eta$=1e-10 and we perform the new-face IR-WRI but using FD as forward engine to avoid inverse crime. The estimated models are shown in Figures \ref{BP_inv_results_IC}a-\ref{BP_inv_results_IC}c. Comparing these results with those shown in Figure \ref{BP_inv_results} highlights the footprint of modelling inaccuracies in the imaging results.  
%

\section{Discussions}
We have proposed a reformulation of frequency-domain FWI with extended search space, which can tackle large-scale 3D problems. 
The search space expansion is generated by computing data-assimilated (DA) wavefields that jointly satisfy the wave equation and the observation equation in a least-squares sense. In the original formulation of IR-WRI, the reconstructed DA wavefields satisfy a normal equation, which is ideally solved with direct methods to process efficiently multiple right-hand sides. 
%
Because of the larger bandwidth of the normal system,  extracting the DA wavefields for 3D large scale problems is not an easy task. In this paper, we use the new face of IR-WRI proposed by \citet{Gholami_2022_EFW}, which reconstructs the DA wavefields by solving several times the classical wave equation instead of solving one time the normal system.
The new algorithm relies on a splitting scheme, which first computes the extended data residuals $\next{\delta \d}^e$ (Eq. \ref{delta_de}) as a deblurred version of the reduced-space data residuals ($\next{\delta \d}^r$), before using them as adjoint sources in Eq. \ref{DAW} to reconstruct the source extensions associated with each physical source. Then, these source extensions are added to the physical sources to compute the DA wavefields with one additional wave-equation solution per source.\\
To mitigate the computational burden generated by this extra number of wave simulations, we use source and receiver sketching methods, which have however some detrimental effects on the reconstructed wavefields and models in terms of convergence speed and accuracy. These effects can be mitigated by filtering out cross-talk noise with sparsity-promoting regularization during model updating. Our numerical results show that the inversion is more sensitive to receiver sketching than source sketching. Therefore, a careful attention should be paid to the adaptive design of receiver sketching with frequency during multiscale inversion.  \\ 
Finally, we propose to use the CBS method as a scalable, limited-memory and accurate forward engine for FWI. CBS can build highly accurate wavefields for 3D large-scale problems when the latter are intractable for direct solvers. 
Put together, the scalability and low-memory demand of the CBS method allow for an efficient parallelization for multi-RHS simulations by distributing RHSs over a large number of processors. This parallelism can be supplemented by performing the fast Fourier transform in parallel if enough resources are available. These are distinct advantages compared to frequency-domain methods based on direct solvers. The time complexity is more difficult to assess as the convergence speed of the CBS method depends on the magnitude of the contrasts. The time complexity of one CBS iteration $n \log(n)$ is roughly the same as that of a FD time-domain (FDTD) method, namely $n$. The number of time steps in FDTD method generally scales to $n^{1/3}$. The number of CBS iterations required to perform a simulation in the true BP salt model with a sufficient accuracy (namely, $\sim$ 2000-3000 according to Figure~\ref{Fig_BP_scalability}a) suggests that the time complexity of CBS and FDTD method is similar in contrasted media, while the CBS method may be one order of magnitude faster than FDTD in low-contrast media such as biological tissues \citep{Osnabrugge_2016_CBS}. 
To support this statement, let's consider the pseudo-speed (the propagated distance during one iteration/time step) of the CBS method and the FDTD method given by $ps_{cbs}=2 / (k_M \nu)$ and $ps_{fdtd}=\pi / 2 k_M \sqrt{D}$, respectively \citep[][ Table 1]{Osnabrugge_2016_CBS}, when $D=3$ in 3D media and the grid interval is one quarter of the minimum wavelength for both methods (this discretization rule is necessary to sample at Nyquist rate an heterogeneity the size of which is half the wavelength, the smallest structure than can be reconstructed by FWI \citep{Virieux_2009_OFW}). For minimum/maximum wavespeeds of 1500~m/s and 4500~m/s and a 3~Hz frequency (grid interval of 125~m), we found $ps_{cbs}$=116.1~m and $ps_{fdtd}$=72.2~m, which have the same order of magnitude.  Let's remind also that the time complexity of FDTD methods and FDFD methods based either on direct or iterative solvers are the same and equal to $\mathcal{O}(n^{2})$ for an acquisition involving $n^{2/3}$ RHSs \citep{Virieux_2009_TLE}.\\
We stress that the new face of IR-WRI can be implemented with any kinds of forward engines. For example, a time-marching method based on finites differences, finite elements or pseudo-spectral approach can be used to perform wavefield simulation in time while the monochromatic wavefields for parameter updating can be generated on the fly in the loop over time steps by discrete Fourier transform or phase sensitive detection \cite{Nihei_2007_FRM,Sirgue_2008_FDW}.

To the best of our knowledge, CBS is only developed for acoustic physics with constant density \citep{Osnabrugge_2016_CBS}. Its extension toward more complicated physics involving heterogeneous density and anisotropy remains a challenge as well as the implmentation of free-surface boundary condition. 
%
%
%
\section{Conclusions}
We revisit iteratively-refined wavefield reconstruction inversion (IR-WRI) such that large-scale 3D seismic imaging problems can be tackled in the frequency domain. The method relies on the classical wave equation rather than the augmented wave equation with the observation equation (resulting in a normal equation), the formed being easier to manage both in the time and frequency domains. However, the new algorithm requires a higher number of wave equation solutions after the reformulation of the normal equation in a form suitable for explicit time-marching methods. We decrease the number of simulations by using sketching methods applied on sources and receivers. Moreover, we computed monochromatic wavefields with the limited-memory discretization-free CBS method to tackle large-scale problems with a high accuracy. Numerical tests show that the new formulation implemented with source and receiver sketching can reconstruct velocity models that are close to those obtained with the original formulation while reducing the number of wave-equation simulations by 70$\%$. Although the CBS method was used in this study, the new formulation can be implemented with any forward engines as finite-difference frequency-domain and time-domain methods or pseudo-spectral methods.
\section*{ACKNOWLEDGMENTS}  
This study was partially funded by the WIND consortium (\textit{https://www.geoazur.fr/WIND}), sponsored by Chevron, Shell and Total. The authors are grateful to the OPAL infrastructure from 
Observatoire de la Côte d'Azur (CRIMSON) for providing resources and support. This work was granted access to the HPC resources of IDRIS under the allocation A0050410596 made by GENCI.

%
\bibliographystyle{cas-model2-names}
%

\newcommand{\SortNoop}[1]{}

\end{document}